\documentclass[nospthms,graybox]{svmult}

\usepackage{mathptmx}        
\usepackage{helvet}          
\usepackage{courier}         
\usepackage{type1cm}         

\usepackage{makeidx}         
\usepackage{graphicx}        
\usepackage{multicol}        
\usepackage[bottom]{footmisc}

\usepackage{amsmath}
\usepackage{amssymb}
\usepackage{amscd}
\usepackage{indentfirst}

\makeindex                   

\usepackage[colorlinks=true,linkcolor=black,citecolor=black,urlcolor=black,pdfpagelabels=false]{hyperref}

\usepackage{mathtools}

\usepackage{amsthm}
\usepackage{thmtools}

\theoremstyle{plain}
  \declaretheorem[numberwithin=section]{theorem}
  \declaretheorem[numberlike=theorem]{entry}
  \declaretheorem[numberlike=theorem]{corollary}

  \declaretheorem[numberlike=theorem]{conjecture}
\theoremstyle{definition}
  \declaretheorem[numberlike=theorem]{example}
  \declaretheorem[numberlike=theorem]{remark}

\numberwithin{equation}{section}

\newcommand{\df}{\dfrac}
\newcommand{\tf}{\tfrac}
\renewcommand{\Re}{\operatorname{Re}}
\renewcommand{\Im}{\operatorname{Im}}
\renewcommand{\a}{\alpha}
\renewcommand{\b}{\beta}

\renewcommand{\(}{\left\(}
\renewcommand{\)}{\right\)}
\renewcommand{\[}{\left\[}
\renewcommand{\]}{\right\]}

\newcommand{\mathd}{\mathrm{d}}

\begin{document}

\title*{Ramanujan's Formula for $\zeta(2n+1)$}
\titlerunning{Ramanujan's Formula for $\zeta(2n+1)$}

\author{Bruce C.~Berndt and Armin Straub}
\authorrunning{B.C.~Berndt and A.~Straub}
\institute{Bruce C.~Berndt \at Department of Mathematics, University of Illinois at Urbana--Champaign,\\ 1409 W. Green St., Urbana, IL 61801, USA\\ \email{berndt@illinois.edu}
\and
Armin Straub \at Department of Mathematics and Statistics, University of South Alabama,\\ 411 University Blvd N, Mobile, AL 36688, USA\\ \email{straub@southalabama.edu}}

\maketitle

\def\theequation{\thesection.\arabic{equation}}

\section{Introduction}

As customary, $\zeta(s)=\sum_{n=1}^{\infty}n^{-s}, \Re s>1$,
denotes the Riemann zeta function.   Let $B_r$,
$r\geq0$, denote the $r$-th Bernoulli number.  When $n$ is a positive
integer, Euler's formula 
\begin{equation}\label{eulerformula}
\zeta(2n)=\df{(-1)^{n-1}B_{2n}}{2(2n)!}(2\pi)^{2n}
\end{equation}
not only provides an elegant formula for evaluating $\zeta(2n)$, but it also tells us of the arithmetical nature of $\zeta(2n)$. In contrast, we know very little about the odd zeta values $\zeta(2n+1)$.  One of the major achievements in number theory in the past half-century is R.~Ap\'ery's proof that $\zeta(3)$ is irrational \cite{apery}, but for $n\geq2$, the arithmetical nature of $\zeta(2n+1)$ remains open.

Ramanujan made many beautiful and elegant discoveries in his short life of 32 years, and one of them that has attracted the attention of several mathematicians over the years is his intriguing formula for $\zeta(2n+1)$.  To be sure, Ramanujan's formula does not possess the elegance of \eqref{eulerformula}, nor does it provide any arithmetical information.  But, one of the goals of this survey is to convince readers that it is indeed a remarkable formula.

\begin{theorem}[Ramanujan's formula for $\zeta(2n+1)$]\label{ie28}
Let $B_r$, $r\geq0$, denote the $r$-th Bernoulli number.
If $ \a $ and $\b$ are positive numbers such that $ \a\b=\pi^2$,
and if $ n $ is a positive integer, then
\begin{align}
&\mspace{-25mu}\a^{-n}\left(\df{1}{2}\zeta(2n+1) +
\sum_{m=1}^\infty \df{1}{m^{2n+1}(e^{2m\a}-1)}\right)\notag\\&\quad-
(-\b)^{-n}\left(\df{1}{2}\zeta(2n+1) +
\sum_{m=1}^\infty \df{1}{m^{2n+1}(e^{2m\b}-1)}\right)\notag\\
&=2^{2n}\sum_{k=0}^{n+1}
(-1)^{k-1}\df{B_{2k}}{(2k)!}\df{B_{2n+2-2k}}{(2n+2-2k)!}\a^{n+1-k}\b^k.
\label{i3.21}
\end{align}
\end{theorem}

Theorem \ref{ie28} appears as Entry 21(i) in Chapter 14 of Ramanujan's second notebook \cite[p.~173]{nb}. It also appears in a formerly unpublished manuscript of Ramanujan that was published in its original handwritten form with his lost notebook \cite[formula (28), pp.~319--320]{lnb}.

The purposes of this paper are to convince readers why \eqref{i3.21} is a fascinating formula, to discuss the history of \eqref{i3.21} and formulas surrounding it, and to discuss the remarkable properties of the polynomials on the right-hand side of \eqref{i3.21}.  We briefly point the readers to analogues and generalizations at the end of our paper.

In Section \ref{sect1}, we discuss Ramanujan's aforementioned unpublished
manuscript and his faulty argument in attempting to prove \eqref{i3.21}.
Companion formulas (both correct and incorrect) with the same parentage are
also examined.  In the following Section \ref{sect2}, we offer an alternative
formulation of \eqref{i3.21} in terms of hyperbolic cotangent sums. 
In Sections~\ref{sec:eisenstein} and \ref{sec:eichler}, we then discuss a
more modern interpretation of Ramanujan's identity.  We introduce Eisenstein
series and their Eichler integrals, and observe that their transformation
properties are encoded in \eqref{i3.21}. 
In particular, this leads us to a vast extension of Ramanujan's identity from
Eisenstein series to general modular forms. 
 
In a different direction, \eqref{i3.21} is a special case of a general transformation formula
for analytic Eisenstein series, or, in another context, a general
transformation formula that greatly generalizes the transformation formula
for the logarithm of the Dedekind eta function.  We show that Euler's famous
formula for $\zeta(2n)$ arises from the same general transformation formula,
and so Ramanujan's formula \eqref{i3.21} is a natural analogue of Euler's
formula.  All of this is discussed in Section \ref{sect3}. 

In Section~\ref{sec:roots}, we discuss some of the remarkable properties of
the polynomials that appear in \eqref{i3.21}.  These polynomials have
received considerable recent attention, with exciting extensions by various
authors to other modular forms.  We sketch recent developments and
indicate opportunities for further research. 
We then provide in Section \ref{sect4} a brief compendium  of proofs of
\eqref{i3.21}, or its equivalent formulation with hyperbolic cotangent sums.

\section{Ramanujan's Unpublished Manuscript}\label{sect1}

The aforementioned handwritten manuscript containing Ramanujan's formula for
$\zeta(2n+1)$ was published for the first time with Ramanujan's lost notebook
\cite[pp.~318--321]{lnb}. 
This partial manuscript was initially examined in detail by the first author
in \cite{jrms}, and by G.E.~Andrews and the first author in their fourth book
on Ramanujan's lost notebook \cite[Chapter 12, pp.~265--284]{geabcbIV}. 

The manuscript's content strongly suggests that it was intended to be a
continuation of Ramanujan's paper \cite{formulae}, \cite[pp.~133--135]{cp}.
It begins with paragraph numbered 18, giving further evidence that it was
intended to be a part of \cite{formulae}.   One might therefore ask why
Ramanujan did not incorporate this partial manuscript in his paper
\cite{formulae}.  As we shall see, one of the primary claims in his
manuscript is false.  Ramanujan's incorrect proof arose from a partial
fraction decomposition, but because he did not have a firm grasp of the
Mittag--Leffler Theorem, he was unable to discern his mis-application of the
theorem.  Most certainly, Ramanujan was fully aware that he had indeed made a
mistake, and consequently he wisely chose not to incorporate the results in
this manuscript in his paper \cite{formulae}.  Ramanujan's mistake arose when
he attempted to find the partial fraction expansion of
$\cot(\sqrt{w\a})\coth(\sqrt{w\b})$.  We now offer  Ramanujan's incorrect
assertion from his formerly unpublished manuscript. 

\begin{entry}[p.~318, formula (21)]\label{ie21}
If $ \a $ and $ \b $ are positive numbers such that $ \a\b=\pi^2$,
then
\begin{equation}\label{i3.1}
\df{1}{2w} + \sum_{m=1}^\infty \left\{\df{m\a\coth(m\a)}{w+m^2\a}+
\df{m\b\coth(m\b)}{w-m^2\b}\right\} =
\df{\pi}{2}\cot(\sqrt{w\a})\coth(\sqrt{w\b}).
\end{equation}
\end{entry}

We do not know if Ramanujan was aware of the Mittag--Leffler Theorem.  Normally, \emph{if}  we could apply this theorem to
$\cot(\sqrt{w\a})\coth(\sqrt{w\b})$, we would let $w\to\infty$ and conclude that the difference between the right- and left-hand sides of \eqref{i3.1} is an entire function that is identically equal to 0.  Of course, in this instance, we cannot make such an argument.
We now offer a corrected version of Entry \ref{ie21}.

\begin{entry}[Corrected Version of \eqref{i3.1}] Under the hypotheses of Entry
  \textup{\ref{ie21}},
\begin{align}\label{i3.10}
\df{\pi}{2}\cot(\sqrt{w\a})\coth(\sqrt{w\b})&= \df{1}{2w} +
\df{1}{2}\log\df{\b}{\a} \\&\quad+
\sum_{m=1}^\infty \left\{\df{m\a\coth(m\a)}{w+m^2\a}+
\df{m\b\coth(m\b)}{w-m^2\b}\right\}.\notag
\end{align}
\end{entry}

Shortly after stating Entry \ref{ie21}, Ramanujan offers the following
corollary. 

\begin{entry}[p.~318, formula (23)]\label{ie23}
If $\a $ and $\b$ are positive numbers such that $ \a\b=\pi^2$,
then
\begin{equation}\label{i3.12}
\a\sum_{m=1}^\infty \df{m}{e^{2m\a}-1}+\b\sum_{m=1}^\infty \df{m}{e^{2m\b}-1}=
\df{\a+\b}{24} -\df{1}{4}.
\end{equation}
\end{entry}

To prove \eqref{i3.12}, it is natural to formally equate coefficients of $ 1/w$ on both
sides of \eqref{i3.1}.  When we do so, we obtain \eqref{i3.12}, but without the term $-\tf14$ on the right-hand side of \eqref{i3.12} \cite[pp.~274--275]{geabcbIV}.  Ramanujan surely must have been perplexed by this, for he had previously proved \eqref{i3.12} by other means. In particular,
Ramanujan offered Entry \ref{ie23} as Corollary (i) in Section 8
of Chapter 14 in his second notebook \cite{nb}, \cite[p.~255]{II}.

Similarly, if one equates constant terms on both sides of \eqref{i3.10}, we obtain the familiar transformation formula for the logarithm of the Dedekind eta function.

\begin{entry}[p.~320, formula (29)]\label{ie29}
If $ \a $ and $\b $ are positive numbers such that $\a\b=\pi^2$,
then
\begin{equation}\label{i3.18}
\sum_{m=1}^\infty \df{1}{m(e^{2m\a}-1)} -
\sum_{m=1}^\infty \df{1}{m(e^{2m\b}-1)}
=\df{1}{4}\log\df{\a}{\b} - \df{\a-\b}{12}.
\end{equation}
\end{entry}

Of course, if we had employed \eqref{i3.1} instead of \eqref{i3.10}, we would
not have obtained the expression $\tf{1}{4}\log\tf{\a}{\b}$ on the right-hand
side of \eqref{i3.18}. 
Entry \ref{ie29} is stated by Ramanujan as Corollary (ii) in
Section 8 of Chapter 14 in his second notebook \cite{nb},
\cite[p.~256]{II} and as Entry 27(iii) in Chapter 16  of his
second notebook \cite{nb}, \cite[p.~43]{III}. 

In contrast to Entries \ref{ie23} and \ref{ie29}, Ramanujan's formula
\eqref{i3.21} can be derived from \eqref{i3.1}, because we  equate
coefficients of $w^n$, $n\geq1$, on both sides, and so the missing terms in
\eqref{i3.1} do not come into consideration.  We now give the argument likely
given by Ramanujan.  (Our exposition is taken from \cite[p.~278]{geabcbIV}.) 

\begin{proof} Return to \eqref{i3.1},  use the equality 
\begin{equation}\label{i3.11}
\coth x = 1+ \df{2}{e^{2x}-1},
\end{equation}
 and    expand the summands into geometric series to arrive at
\begin{align}\label{i3.22}
\mspace{-25mu}\df{\pi}{2}\cot(\sqrt{w\a})\coth(\sqrt{w\b})&= \df{1}{2w}
+\sum_{m=1}^\infty \left\{\df{1}{m}\sum_{k=0}^\infty \left(-\df{w}{m^2\a}\right)^k
\left(1+\df{2}{e^{2m\a}-1}\right)\right.\notag\\&\quad\left.-
\df{1}{m}\sum_{k=0}^\infty \left(\df{w}{m^2\b}\right)^k
\left(1+\df{2}{e^{2m\b}-1}\right)\right\}.
\end{align}
Following Ramanujan, we equate coefficients of $w^n$, $n
\geq 1$, on both sides of \eqref{i3.22}. On the right side, the
coefficient of $w^n$ equals
\begin{align}\label{i3.23}
\mspace{-25mu}(-\a)^{-n}\zeta(2n+1)&+2(-\a)^{-n}\sum_{m=1}^\infty \df{1}{m^{2n+1}(e^{2m\a}-1)}
\notag \\
-\b^{-n}\zeta(2n+1)&+2\b^{-n}\sum_{m=1}^\infty \df{1}{m^{2n+1}(e^{2m\b}-1)}.
\end{align}
Using the Laurent expansions for $\cot z$ and $\coth z $ about $
z=0$, we find that on the left side of \eqref{i3.22},
\begin{align}\label{i3.24}
\df{\pi}{2}\cot(\sqrt{w\a})\coth(\sqrt{w\b})&=\df{\pi}{2}
\sum_{k=0}^\infty \df{(-1)^k2^{2k}B_{2k}}{(2k)!}(w\a)^{k-1/2}\notag\\
&\quad\times\sum_{j=0}^\infty \df{2^{2j}B_{2j}}{(2j)!}(w\b)^{j-1/2}.
\end{align}
The coefficient of $ w^n$ in \eqref{i3.24} is easily seen to be
equal to
\begin{equation}\label{i3.25}
2^{2n+1}\sum_{k=0}^{n+1}(-1)^k\df{B_{2k}}{(2k)!}\df{B_{2n+2-2k}}
{(2n+2-2k)!}\a^k\b^{n+1-k},
\end{equation}
where we used the equality $\a\b=\pi^2$. Now equate the
expressions in \eqref{i3.23} and \eqref{i3.25}, then multiply both
sides by $(-1)^n\tf12$, and lastly replace $k $ by $ n+1-k$ in
the finite sum.
The identity \eqref{i3.21} immediately follows. 
\end{proof}

As mentioned in the Introduction, Ramanujan also recorded \eqref{i3.21} as  Entry 21(i) of Chapter 14 in his second notebook \cite[p.~173]{nb}, \cite[p.~271]{II}.  Prior to that on page 171, he offered the partial fraction decomposition
\begin{equation}\label{pfd}
\pi^2xy\cot(\pi x)\coth(\pi y)=1+2\pi xy\sum_{n=1}^\infty \df{n\,\coth(\pi nx/y)}{n^2+y^2}-
2\pi xy\sum_{n=1}^\infty \df{n\,\coth(\pi ny/x)}{n^2-x^2},
\end{equation}
which should be compared with \eqref{i3.1}.  The two infinite series on the right side of \eqref{pfd} diverge individually, but when combined together into one series, the series converges.  Thus, most likely, Ramanujan's derivation of \eqref{i3.21} when he recorded Entry 21(i) in his second notebook was similar to the argument that he used in his partial manuscript.

R.~Sitaramachandrarao \cite{sita} modified Ramanujan's approach via partial fractions to avoid the pitfalls experienced by Ramanujan. Sitaramachandrarao established the partial fraction decomposition
\begin{align}\label{i3.8}
\pi^2xy\cot(\pi x)\coth(\pi y) &= 1+\df{\pi^2}{3}(y^2-x^2) \\
&\quad-2\pi xy\sum_{m=1}^\infty \left(\df{y^2\coth(\pi
mx/y)}{m(m^2+y^2)} + \df{x^2\coth(\pi my/x)}{m(m^2-x^2)}\right).
\notag
\end{align}
Using the elementary identities
\begin{equation*}
\df{y^2}{m(m^2+y^2)} = -\df{m}{m^2+y^2}+\df{1}{m}
\end{equation*}
and
\begin{equation*}
\df{x^2}{m(m^2-x^2)} = \df{m}{m^2-x^2}-\df{1}{m},
\end{equation*}
and then employing \eqref{i3.11}, he showed that
{\allowdisplaybreaks\begin{align}\label{i3.9}
\pi^2xy\cot(\pi x)\coth(\pi y)
&= 1+\df{\pi^2}{3}(y^2-x^2) \\
&\quad+2\pi xy\sum_{m=1}^\infty \left(\df{m\coth(\pi mx/y)}{m^2+y^2}
-
\df{m\coth(\pi my/x)}{m^2-x^2}\right) \notag\\
&\quad-4\pi xy\sum_{m=1}^\infty \df{1}{m}\left(\df{1}{e^{2\pi
mx/y}-1} - \df{1}{e^{2\pi my/x}-1}\right).\notag
\end{align}}%
Setting $\pi x=\sqrt{w\a}$ and $\pi y=\sqrt{w\b}$ above and invoking the
transformation formula for the logarithm of the  Dedekind eta function from
Entry \ref{ie29}, we can readily deduce \eqref{i3.10}.  For more details, see
\cite[pp.~272--273]{geabcbIV}.

The first published proof of \eqref{i3.21} is due to S.L.~Malurkar
\cite{malurkar}
in 1925--1926.  Almost
certainly, he was unaware that \eqref{i3.21} can be found in Ramanujan's
notebooks \cite{nb}. 

If we set $\a=\b=\pi$ and replace $n$ by $2n+1$ in \eqref{i3.21}, we deduce
that, for $n\geq0$, 
\allowdisplaybreaks{\begin{align}
\zeta(4n+3)=2^{4n+2}\pi^{4n+3}\sum_{k=0}^{2n+2}(-1)^{k+1}\df{B_{2k}}{(2k)!}\df{B_{4n+4-2k}}{(4n+4-2k)!}
-2\sum_{k=1}^\infty \df{k^{-4n-3}}{e^{2\pi k}-1}.
\label{i3.21*}
\end{align}}%
This special case is actually due to M.~Lerch \cite{lerch} in 1901.  The
identity \eqref{i3.21*} is a remarkable identity, for it shows that
$\zeta(4n+3)$ is equal to a rational multiple of $\pi^{4n+3}$ plus a very
rapidly convergent series.  Therefore, $\zeta(4n+3)$ is ``almost'' a rational
multiple of $\pi^{4n+3}$.

\section{An Alternative Formulation of \eqref{i3.21} in Terms of Hyperbolic Cotangent Sums}\label{sect2}

If we use the elementary identity \eqref{i3.11}, we can transform Ramanujan's identify  \eqref{i3.21} for $\zeta(2n+1)$  into an identity for hyperbolic cotangent sums, namely,
\begin{multline}\label{cothangent}
\a^{-n}\sum_{m=1}^\infty \df{\coth(\a m)}{m^{2n+1}}=(-\b)^{-n}\sum_{m=1}^\infty \df{\coth(\b m)}{m^{2n+1}}\\
-2^{2n+1}\sum_{k=0}^{n+1}(-1)^k\df{B_{2k}}{(2k)!}\df{B_{2n+2-2k}}{(2n+2-2k)!}\a^{n+1-k}\b^{k},
\end{multline}
where, as before, $n$ is a positive integer and $\a\b=\pi^2$.
To the best of our knowledge, the first recorded proof of Ramanujan's formula \eqref{i3.21} in the form \eqref{cothangent} was by T.S.~Nanjundiah \cite{nanjundiah} in 1951.
If we replace $n$ by $2n+1$ and set $\a=\b=\pi$, then \eqref{cothangent} reduces to
\begin{equation}\label{cot}
\sum_{m=1}^\infty \df{\coth(\pi m)}{m^{4n+3}}=2^{4n+2}{\pi^{4n+3}}\sum_{k=0}^{2n+2}(-1)^{k+1}\df{B_{2k}}{(2k)!}\df{B_{4n+4-2k}}{(4n+4-2k)!}.
\end{equation}
This variation of \eqref{i3.21*} was also first established by Lerch \cite{lerch}.  Later proofs were given by G.N.~Watson \cite{watsonII}, H.F.~Sandham \cite{sandham}, J.R.~Smart \cite{smart}, F.P.~Sayer \cite{sayer}, Sitaramachandrarao \cite{sita}, and the first author \cite{dedekindhardy}, \cite{rocky}.
The special cases $n=0$ and $n=1$ are Entries 25(i), (ii) in Chapter 14 of
Ramanujan's second notebook \cite[p.~176]{nb}, \cite[p.~293]{II}.  He
communicated this second special case in his first letter to Hardy
\cite[p.~xxvi]{cp}.

 Significantly generalizing an idea of C.L.~Siegel in
deriving the
transformation formula for the Dedekind eta function \cite{siegel},
S.~Kongsiriwong \cite{kongsiriwong} not only established \eqref{cothangent}, but
he derived several beautiful generalizations and analogues  of
\eqref{cothangent}.

Deriving a very general transformation formula for Barnes' multiple zeta
function, Y.~Komori, K.~Matsumoto, and H.~Tsumura \cite{kmt} established not only
\eqref{cothangent}, but also a plethora of further identities and summations
in closed form for series involving hyperbolic trigonometric functions.

\section{Eisenstein Series}\label{sec:eisenstein}

A more modern interpretation of Ramanujan's remarkable formula, which is
discussed for instance in \cite{gunmurtyrath}, is that \eqref{i3.21}
encodes the fundamental transformation properties of Eisenstein series of
level~$1$ and their Eichler integrals. The goal of this section and the next
is to explain this connection and to introduce these modular objects starting
with Ramanujan's formula. As in \cite{gross1b} and
\cite{gunmurtyrath}, set
\begin{equation}
  F_a (z) = \sum_{n = 1}^{\infty} \sigma_{- a} (n) e^{2 \pi i n z},
  \quad \sigma_k (n) = \sum_{d|n} d^k . \label{eq:Fa}
\end{equation}
Observe that, with $q = e^{2 \pi i z}$, we can express $F_a (z)$ as a
Lambert series
\begin{equation*}
  F_a (z) = \sum_{n = 1}^{\infty} \left(\sum_{d|n} d^{- a} \right) q^n =
   \sum_{d = 1}^{\infty} \sum_{m = 1}^{\infty} d^{- a} q^{d m} = \sum_{n =
   1}^{\infty} \frac{n^{- a} q^n}{1 - q^n} = \sum_{n = 1}^{\infty} \frac{n^{-
   a}}{e^{- 2 \pi i n z} - 1}
\end{equation*}
in the form appearing in Ramanujan's formula \eqref{i3.21}. Indeed, if we
let $z = \alpha i / \pi$, then Ramanujan's formula \eqref{i3.21}
translates into
\begin{eqnarray}
  &  & \left\{ \frac{\zeta (2 m + 1)}{2} + F_{2 m + 1} (z) \right\} =
  z^{2 m} \left\{ \frac{\zeta (2 m + 1)}{2} + F_{2 m + 1} \left(-
  \frac{1}{z} \right) \right\} \nonumber\\
  &  & + \frac{(2 \pi i)^{2 m + 1}}{2 z} \sum_{n = 0}^{m + 1} \frac{B_{2
  n}}{(2 n) !}  \frac{B_{2 m - 2 n + 2}}{(2 m - 2 n + 2) !} z^{2 n} . 
  \label{eq:rama:F}
\end{eqnarray}
This generalization to values $z$ in the upper half-plane $\mathcal{H}= \{
z \in \mathbb{C}: \Im (z) > 0 \}$ was derived by E.~Grosswald in
\cite{gross1}.

Ramanujan's formula becomes particularly simple if the integer $m$ satisfies
$m < - 1$. In that case, with $k = - 2 m$, equation \eqref{eq:rama:F} can be
written as
\begin{equation}
  E_k (z) = z^{- k} E_k (- 1 / z), \label{eq:E:S}
\end{equation}
with
\begin{equation}
  E_k (z) = 1 + \frac{2}{\zeta (1 - k)} F_{1 - k} (z) = 1 +
  \frac{2}{\zeta (1 - k)} \sum_{n = 1}^{\infty} \sigma_{k - 1} (n) q^n .
  \label{eq:E}
\end{equation}
The series $E_k$, for even $k > 2$, are known as the {\emph{normalized
Eisenstein series}} of weight $k$. They are fundamental instances of modular
forms. A {\emph{modular form}} of weight $k$ and level~$1$ is a function $G$
on the upper half-plane $\mathcal{H}$, which is holomorphic on $\mathcal{H}$
(and as $z \rightarrow i \infty$), and transforms according to
\begin{equation}
  (c z + d)^{- k} G (V z) = G (z) \label{eq:MF:g}
\end{equation}
for each matrix
\begin{equation*}
  V = \begin{bmatrix}
     a & b\\
     c & d
   \end{bmatrix} \in \operatorname{SL}_2 (\mathbb{Z}) .
\end{equation*}
Here, as usual, $\operatorname{SL}_2 (\mathbb{Z})$ is the modular group of integer
matrices $V$ with determinant~$1$, and these act on $\mathcal{H}$ by
fractional linear transformations as
\begin{equation}
  V z = \frac{a z + b}{c z + d} . \label{eq:V}
\end{equation}
Modular forms of higher level are similarly defined and transform only with
respect to certain subgroups of $\operatorname{SL}_2 (\mathbb{Z})$. Equation
\eqref{eq:E:S}, together with the trivial identity $E_k (z + 1) = E_k (z)$, establishes that $E_k (z)$ satisfies the modular transformation
\eqref{eq:MF:g} for the matrices
\begin{equation*}
  S = \begin{bmatrix}
     0 & - 1\\
     1 & 0
   \end{bmatrix}, \quad T = \begin{bmatrix}
     1 & 1\\
     0 & 1
   \end{bmatrix} .
\end{equation*}
It is well-known that $S$ and $T$ generate the modular group $\operatorname{SL}_2 (\mathbb{Z})$, from which it follows that $E_k (z)$ is invariant under the
action of any element in $\operatorname{SL}_2 (\mathbb{Z})$. In other words, $E_k (z)$ satisfies \eqref{eq:MF:g} for all matrices in $\operatorname{SL}_2 (\mathbb{Z})$.

To summarize our discussion so far, the cases $m < - 1$ of Ramanujan's formula
\eqref{i3.21} express the fact that, for even $k > 2$, the $q$-series $E_k (z)$, defined in \eqref{eq:E} as the generating function of sums of powers
of divisors, transforms like a modular form of weight $k$. Similarly, the case
$m = - 1$ encodes the fact that
\begin{equation*}
  E_2 (z) = 1 - 24 \sum_{n = 1}^{\infty} \sigma_1 (n) q^n,
\end{equation*}
as defined in \eqref{eq:E}, transforms according to
\begin{equation*}
  E_2 (z) = z^{- 2} E_2 (- 1 / z) - \frac{6}{\pi i z} .
\end{equation*}
This is an immediate consequence of Ramanujan's formula in the form
\eqref{eq:rama:F}. The weight $2$ Eisenstein series $E_2 (z)$ is the
fundamental instance of a {\emph{quasimodular form}}.

In the next section, we discuss the case $m > 0$ of Ramanujan's formula
\eqref{i3.21} and describe its connection with Eichler integrals.

\section{Eichler Integrals and Period Polynomials}\label{sec:eichler}

If $f$ is a modular form of weight $k$ and level~1, then
\begin{equation*}
  z^{- k} f \left(- \frac{1}{z} \right) - f (z) = 0.
\end{equation*}
The derivatives of $f$, however, do not in general transform in a modular
fashion by themselves. For instance, taking the derivative of this modular
relation and rearranging, we find that
\begin{equation*}
  z^{- k - 2} f' \left(- \frac{1}{z} \right) - f' (z) =
   \frac{k}{z} f (z) .
\end{equation*}
This is a modular transformation law only if $k = 0$, in which case the
derivative $f'$ (or, the $- 1$st integral of $f$) transforms like a modular
form of weight $2$. Remarkably, it turns out that any $(k - 1)$st integral of
$f$ does satisfy a nearly modular transformation rule of weight $k - 2 (k -
1) = 2 - k$. Indeed, if $F$ is a $(k - 1)$st primitive of $f$, then
\begin{equation}
  z^{k - 2} F \left(- \frac{1}{z} \right) - F (z)
  \label{eq:eichler:poly:intro}
\end{equation}
is a polynomial of degree at most $k - 2$. The function $F$ is called an
{\emph{Eichler integral}} of $f$ and the polynomials are referred to as
{\emph{period polynomials}}.

Let, as usual,
\begin{equation*}
  D = \frac{1}{2 \pi i}  \frac{\mathd}{\mathd z} = q \frac{\mathd}{\mathd
   q} .
\end{equation*}
\begin{remark}
  That \eqref{eq:eichler:poly:intro} is a polynomial can be seen as a
  consequence of {\emph{Bol's identity}} \cite{bol}, which states that,
  for all sufficiently differentiable $F$ and $V = \begin{bsmallmatrix}
    a & b\\
    c & d
  \end{bsmallmatrix} \in \operatorname{SL}_2 (\mathbb{R})$,
  \begin{equation}
    \frac{(D^{k - 1} F) (V z)}{(c z + d)^k} = D^{k - 1} [ (c z +
    d)^{k - 2} F (V z)] . \label{eq:bol}
  \end{equation}
  Here, $V z$ is as in \eqref{eq:V}. If $F$ is an Eichler integral of a
  modular form $f = D^{k - 1} F$ of weight $k$, then
  \begin{equation*}
    D^{k - 1} [ (c z + d)^{k - 2} F (V z) - F (z)] = \frac{(D^{k -
     1} F) (V z)}{(c z + d)^k} - (D^{k - 1} F) (z) = 0,
  \end{equation*}
  for all $V \in \operatorname{SL}_2 (\mathbb{Z})$. This shows that $(c z + d)^{k
  - 2} F (V z) - F (z)$ is a polynomial of degree at most $k - 2$. A
  delightful explanation of Bol's identity \eqref{eq:bol} in terms of Maass
  raising operators is given in \cite[IV.~2]{lz}.
\end{remark}

Consider, for $m > 0$, the series
\begin{equation*}
  F_{2 m + 1} (z) = \sum_{n = 1}^{\infty} \sigma_{- 2 m - 1} (n) q^n,
\end{equation*}
which was introduced in \eqref{eq:Fa} and which is featured in Ramanujan's
formula \eqref{eq:rama:F}. Observe that
\begin{equation*}
  D^{2 m + 1} F_{2 m + 1} (z) = \sum_{n = 1}^{\infty} \left(\sum_{d|n}
   d^{- 2 m - 1} \right) n^{2 m + 1} q^n = \sum_{n = 1}^{\infty} \sigma_{2 m +
   1} (n) q^n,
\end{equation*}
and recall from \eqref{eq:E} that this is, up to the missing constant term, an
Eisenstein series of weight $2 m + 2$. It thus becomes clear that Ramanujan's
formula \eqref{eq:rama:F}, in the case $m > 0$, encodes the fundamental
transformation law of the Eichler integral corresponding to the weight $2 m +
2$ Eisenstein series.

\begin{remark}
  \label{rk:pp}That the right-hand side of \eqref{eq:rama:F} is only almost a
  polynomial is due to the fact that its left-hand side needs to be adjusted
  for the constant term of the underlying Eisenstein series. To be precise,
  integrating the weight $k = 2 m + 2$ Eisenstein series \eqref{eq:E},
  slightly scaled, $k - 1$ times, we find that
  \begin{equation}
    G_{2 m + 1} (z) := \frac{\zeta (- 1 - 2 m)}{2}  \frac{(2 \pi i
    z)^{2 m + 1}}{(2 m + 1) !} + F_{2 m + 1} (z) \label{eq:ei:G}
  \end{equation}
  is an associated Eichler integral of weight $- 2 m$. Keeping in mind the
  evaluation $\zeta (- 1 - 2 m) = - B_{2 m + 2} / (2 m + 2)$, Ramanujan's
  formula in the form \eqref{eq:rama:F} therefore can be restated as
  \begin{eqnarray}
    z^{2 m} G_{2 m + 1} \left(- \frac{1}{z} \right) - G_{2 m + 1} (z) & = & \frac{\zeta (2 m + 1)}{2} (1 - z^{2 m})  \label{eq:pp:G}\\
    &  & - \frac{(2 \pi i)^{2 m + 1}}{2} \sum_{n = 1}^m \frac{B_{2 n}}{(2 n)
    !}  \frac{B_{2 m - 2 n + 2}}{(2 m - 2 n + 2) !} z^{2 n - 1}, \nonumber
  \end{eqnarray}
  where the right-hand side is a polynomial of degree $k - 2 = 2 m$. Compare,
  for instance, \cite[eq.~(11)]{zagier91}.
\end{remark}

An interesting and crucial property of period polynomials is that their
coefficients encode the critical $L$-values of the corresponding modular form.
Indeed, consider a modular form
\begin{equation*}
  f (z) = \sum_{n = 0}^{\infty} a (n) q^n
\end{equation*}
of weight $k$. For simplicity, we assume that $f$ is modular on the full
modular group $\operatorname{SL}_2 (\mathbb{Z})$ (for the case of higher level, we
refer to \cite{pp}). Let
\begin{equation*}
  f^{\ast} (z) := \frac{a (0) z^{k - 1}}{(k - 1) !} +
   \frac{1}{(2 \pi i)^{k - 1}} \sum_{n = 1}^{\infty} \frac{a (n)}{n^{k - 1}}
   q^n
\end{equation*}
be an Eichler integral of $f$. Then, a special case of a result by M.J.~Razar
\cite{razar} and A.~Weil \cite{weil} (see also
\cite[(8)]{gunmurtyrath}) shows the following.

\begin{theorem}
  If $f$ and $f^{\ast}$ are as described above, then
  \begin{equation}
    z^{k - 2} f^{\ast} \left(- \frac{1}{z} \right) - f^{\ast} (z) =
    - \sum_{j = 0}^{k - 2} \frac{L (k - 1 - j, f)}{j! (2 \pi i)^{k - 1 - j}}
    z^j, \label{eq:rw}
  \end{equation}
  where, as usual, $L (s, f) = \sum_{n = 1}^{\infty} a_n n^{- s}$.
\end{theorem}

\begin{example}
  For instance, if
  \begin{equation*}
    f (z) = \frac{1}{2} (2 \pi i)^{k - 1} \zeta (1 - k) E_k (z) = (2 \pi i)^{k - 1} \left[ \frac{1}{2} \zeta (1 - k) + \sum_{n =
     1}^{\infty} \sigma_{k - 1} (n) q^n \right],
  \end{equation*}
  where $E_k$ is the Eisenstein series defined in \eqref{eq:E}, then
  \begin{equation*}
    f^{\ast} (z) = \frac{1}{2} (2 \pi i)^{k - 1} \zeta (1 - k)
     \frac{z^{k - 1}}{(k - 1) !} + F_{k - 1} (z)
  \end{equation*}
  equals the function $G_{2 m + 1} (z)$, where $k = 2 m + 2$, used earlier
  in \eqref{eq:ei:G}. Observe that the $L$-series of the Eisenstein series $f
  (z)$ is
  \begin{equation*}
    L (s, f) = (2 \pi i)^{k - 1} \sum_{n = 1}^{\infty} \frac{\sigma_{k - 1}
     (n)}{n^s} = (2 \pi i)^{k - 1} \zeta (s) \zeta (1 - k + s) .
  \end{equation*}
  The result \eqref{eq:rw} of Razar and Weil therefore implies that
  \begin{equation}
    z^{k - 2} f^{\ast} \left(- \frac{1}{z} \right) - f^{\ast} (z) =
    - \sum_{j = 0}^{k - 2} \frac{(2 \pi i)^j \zeta (k - 1 - j) \zeta (-
    j)}{j!} z^j, \label{eq:rama:rw}
  \end{equation}
  and the right-hand side of \eqref{eq:rama:rw} indeed equals the right-hand
  side of \eqref{eq:pp:G}, which we obtained as a reformulation of Ramanujan's
  identity. In particular, we see that \eqref{eq:rw} is a vast generalization
  of Ramanujan's identity from Eisenstein series to other modular forms.
\end{example}

\section{Ramanujan's Formula for $\zeta(2n+1)$ and Euler's Formula for $\zeta(2n)$ are Consequences of the Same General Theorem}\label{sect3}

Ramanujan's formula \eqref{i3.21} for $\zeta(2n+1)$ is, in fact, a special
instance of a general transformation formula for generalized analytic
Eisenstein series or, in another formulation, for a vast generalization of
the logarithm of the Dedekind eta function.  We relate one such
generalization due to the first author \cite{bcbtrans} and developed in
\cite{rocky}, where a multitude of special cases are derived.  Throughout,
let $z$ be in the upper half-plane $\mathcal{H}$,
and let $V = \begin{bsmallmatrix} a & b\\ c & d \end{bsmallmatrix}$ be an element of $\operatorname{SL}_2(\mathbb{Z})$.
That is, $a,b,c,d\in\mathbb{Z}$ and $ad-bc=1$.
Further, let $Vz=(az+b)/(cz+d)$, as in \eqref{eq:V}.
In the sequel, we assume that $c>0$. Let $r_1$ and $r_2$ be real numbers, and let $m$ be an integer.
Define
\begin{equation}\label{2.1}
 A(z,-m,r_1,r_2):=\sum_{n>-r_1}\sum_{k=1}^\infty k^{-m-1}e^{2\pi
   ik(nz+r_1z+r_2)}
\end{equation}
and
\begin{equation}\label{2.2}
 H(z,-m,r_1,r_2):= A(z,-m,r_1,r_2)+(-1)^m A(z,-m,-r_1,-r_2).
\end{equation}
We define the Hurwitz zeta function $\zeta(s,a)$, for any real number $a$, by
$$\zeta(s,a):=\sum_{n>-a}(n+a)^{-s}, \quad \Re s>1.$$
For real numbers $a$, the characteristic function is denoted by
$\lambda(a)$. Define $R_1=ar_1+cr_2$ and $R_2=br_1+dr_2$, where $a,b,c,d$ are
as above.  Define
\begin{align}\label{g}
  g(z,-m,r_1,r_2)&:=\lim_{s\to-m} \frac{\Gamma(s)}{(2\pi i)^s}
  \biggl\{-\lambda(r_1)(cz+d)^{-s}(e^{\pi is}\zeta(s,r_2)+\zeta(s,-r_2))\notag\\
&\quad +\lambda(R_1)(\zeta(s,R_2)+e^{-\pi is}\zeta(s,-R_2))\biggr\}.
\end{align}
As customary, $B_n(x)$, $n\geq0$, denotes the $n$-th Bernoulli polynomial,
and $\{x\}$ denotes the fractional part of $x$.
Let $\rho:=\{R_2\}c-\{R_1\}d$. Lastly, define
\begin{equation*}
h(z,-m,r_1,r_2):=\sum_{j=1}^c\sum_{k=0}^{m+2}
\frac{(-1)^{k-1}(cz+d)^{k-1}}{k!(m+2-k)!}
B_k\left(\df{j-\{R_1\}}{c}\right)B_{m+2-k}\left(\df{jd+\rho}{c}\right),
\end{equation*}
where it is understood that if $m+2<0$, then $h(z,-m,r_1,r_2)\equiv0$.

We are now ready to state our general transformation formula
\cite[p.~498]{bcbtrans}, \cite[p.~150]{rocky}.

\begin{theorem}\label{thmH} If $z\in\mathcal{H}$ and $m$ is any integer, then
\begin{align}
(cz+d)^mH(Vz,-m,r_1,r_2)&=H(z,-m,R_1,R_2)\notag\\&\quad+g(z,-m,r_1,r_2)+(2\pi
i)^{m+1}h(z,-m,r_1,r_2). \label{2.6}
\end{align}
\end{theorem}

We now specialize Theorem \ref{thmH} by supposing that $r_1=r_2=0$ and that
$Vz=-1/z$, so that $c=1$ and $d=0$.
Note that, from \eqref{2.1} and \eqref{2.2},
\begin{equation}\label{2.9}
H(z,-m,0,0)=(1+(-1)^m)\sum_{k=1}^\infty \df{k^{-m-1}}{e^{-2\pi ikz}-1}.
\end{equation}

\begin{theorem}\label{thmHH}
If $z\in\mathcal{H}$ and $m$ is any integer, then
\begin{align}\label{2.10}
&z^m(1+(-1)^m)\sum_{k=1}^\infty \df{k^{-m-1}}{e^{2\pi ik/z}-1}=
(1+(-1)^m)\sum_{k=1}^\infty \df{k^{-m-1}}{e^{-2\pi ikz}-1}+g(z,-m)\notag\\
&\qquad+(2\pi
i)^{m+1}\sum_{k=0}^{m+2}\df{B_k(1)}{k!}\df{B_{m+2-k}}{(m+2-k)!}(-z)^{k-1},
\end{align}
where
\begin{equation}
g(z,-m)=\begin{cases}\pi i-\log z, \quad &\text{if $m=0$},\\
\{1-(-z)^m\}\zeta(m+1),&\text{if $m \neq0$}.\end{cases}\label{2.7,2.8}
\end{equation}
\end{theorem}

\begin{corollary}[Euler's Formula for $\zeta(2n)$]\label{e} For each positive
  integer   $n$,
\begin{equation}\label{eulerzeta}
\zeta(2n)=\df{(-1)^{n-1}B_{2n}}{2(2n)!}(2\pi)^{2n}.
\end{equation}
\end{corollary}

\begin{proof} Put $m=2n-1$ in \eqref{2.10}.  Trivially, by \eqref{2.9},
  $H(z,-2m+1,0,0)=0$. Using
  \eqref{2.7,2.8}, we see that \eqref{2.10} reduces to
\begin{equation}\label{euler}
(1+z^{2n-1})\zeta(2n)=
\df{(2\pi)^{2n}(-1)^{n-1}}{(2n)!}\{B_1(1)B_{2n}-B_{2n}B_1z^{2n-1}\},
\end{equation}
where we have used the values $B_k(1)=B_k$, $k\geq2$, and $B_{2k+1}=0$, $k\geq1$.  Since $B_1(1)=\tf12$
and $B_1=-\tf12$, Euler's formula \eqref{eulerzeta} follows immediately from
\eqref{euler}.
\end{proof}

\begin{corollary}[Ramanujan's Formula for $\zeta(2n+1)$]\label{r} Let $\a$
  and $\b$   denote positive numbers such that $\a\b=\pi^2$. Then, for each
  positive   integer $n$,
\begin{align}\label{2.11}
&\a^{-n}\left\{\df12\zeta(2n+1)+
\sum_{k=1}^\infty \df{k^{-2n-1}}{e^{2\a k}-1}\right\}
=(-\b)^{-n}\left\{\df12\zeta(2n+1)+
\sum_{k=1}^\infty \df{k^{-2n-1}}{e^{2\b k}-1}\right\}\notag\\
&\qquad-2^{2n}\sum_{k=0}^{n+1}(-1)^k\df{B_{2k}}{(2k)!}
\df{B_{2n+2-2k}}{(2n+2-2k)!}\a^{n+1-k}\b^k.
\end{align}
\end{corollary}

\begin{proof}
Set $m=2n$ in \eqref{2.10}, and let $z=i\pi/\a$.  Recall that $\pi^2/\a=\b$.  If we multiply both
sides by $\tf12(-\b)^{-n}$, we obtain \eqref{2.11}.
\end{proof}

We see from Corollaries \ref{e} and \ref{r} that Euler's formula for
$\zeta(2n)$ and Ramanujan's formula for $\zeta(2n+1)$ are natural companions,
because both are special instances of the same transformation
formula.

We also emphasize that in the foregoing proof, we assumed that $n$
was a \emph{positive} integer.  Suppose that we had taken $m=-2n$, where $n$
is a positive integer.  Furthermore, if $n>1$, then the sum
of terms involving Bernoulli
numbers is empty.  Recall that
\begin{equation}\label{zetaneg}
\zeta(1-2n)=-\df{B_{2n}}{2n}, \quad n\geq 1.
\end{equation}
We then obtain the following corollary.

\begin{corollary}\label{m} Let $\a$ and $\b$ be positive numbers such that
  $\a\b=\pi^2$. Then, for any integer $n>1$,
\begin{equation}\label{2.16}
\a^n\sum_{k=1}^\infty \df{k^{2n-1}}{e^{2\a k}-1}-
(-\b)^n\sum_{k=1}^\infty \df{k^{2n-1}}{e^{2\b k}-1}
=\{\a^n-(-\b)^n\}\df{B_{2n}}{4n}.
\end{equation}
\end{corollary}

Corollary \ref{m} is identical to Entry 13 in Chapter 14 of Ramanujan's
second notebook \cite{nb}, \cite[p.~261]{II}. It can also be found in his
paper \cite[p.~269]{trigsums}, \cite[p.~190]{cp} stated without
proof. Corollary \ref{m} is also formula (25) in Ramanujan's unpublished
manuscript \cite[p.~319]{lnb}, \cite[p.~276]{geabcbIV}. 

 Of course, by \eqref{zetaneg}, we could regard \eqref{2.16} as a formula for
 $\zeta(1-2n)$, and 
 so we would have a third  companion to the formulas of Euler and
Ramanujan. 
 The first proof of \eqref{2.16} known to the authors is by
M.B.~Rao and M.V.~Ayyar \cite{ayyar} in 1923, with Malurkar \cite{malurkar}
giving another proof shortly thereafter in 1925.
 If we replace $n$ by $2n+1$ in
\eqref{2.16}, where $n$ is a positive integer, and set $\a=\b=\pi$, we obtain
the special case 
\begin{equation}
\sum_{k=1}^\infty \df{k^{4n+1}}{e^{2\pi k}-1}=\df{B_{4n+2}}{4(2n+1)},\label{glaisher}
\end{equation}
which is due much earlier to J.W.L.~Glaisher \cite{glaisher} in 1889.  The formula \eqref{glaisher} can also be found in Section 13 of Chapter 14 in Ramanujan's second notebook \cite{nb}, \cite[p.~262]{II}.

There is yet a fourth companion.  If we set $m=-2$ in \eqref{2.10}, 
proceed as in the proof of Corollary \ref{r}, and use \eqref{zetaneg}, we
deduce the following 
corollary, which we have previously recorded as Entry \ref{ie23}, and
which can be thought of as ``a formula for $\zeta(-1)$.'' 

\begin{corollary}\label{s}
Let $\a$ and $\b$ be positive numbers such that $\a\b=\pi^2$. Then
\begin{equation}\label{2.20}
\a\sum_{k=1}^\infty \df{k}{e^{2\a k}-1}+\b\sum_{k=1}^\infty \df{k}{e^{2\b k}-1}=
\df{\a+\b}{24}-\df{1}{4}.
\end{equation}
\end{corollary}

If $\a=\b=\pi$, \eqref{2.20} reduces to
\begin{equation}\label{2.21}
\sum_{k=1}^\infty \df{k}{e^{2\pi k}-1}=\df{1}{24}-\df{1}{8\pi}.
\end{equation}

Both the special case \eqref{2.21} and the more general identity \eqref{2.20}
can be found in Ramanujan's notebooks \cite[vol.~1, 
p.~257, no.~9; vol.~2, p.~170, Cor.~1]{nb}, \cite[pp.~255--256]{II}.
However, in 1877, O.~Schl\"omilch 
\cite{sch1}, \cite[p.~157]{sch2} apparently gave the first proof of both
\eqref{2.21} and \eqref{2.20}.

In conclusion of Section \ref{sect3}, we point out that several authors have
proved general transformation formulas from which Ramanujan's formula
\eqref{2.11} for 
$\zeta(2n+1)$  can be deduced as a special case.  However, in
most cases, Ramanujan's formula was not explicitly recorded by the authors.
General transformation formulas have been proved by A.P.~Guinand
\cite{guinand1}, \cite{guinand2}, K.~Chandrasekharan and R.~Narasimhan
\cite{cn}, T.M.~Apostol \cite{apostol1}, \cite{apostol2}, M.~Mikol\'as
\cite{mikolas}, K.~Iseki \cite{iseki}, R.~Bodendiek \cite{bodendiek},
Bodendiek and U.~Halbritter \cite{halbritter},  H.-J.~Glaeske
\cite{glaeske1}, \cite{glaeske2}, D.M.~Bradley \cite{bradley},  Smart
\cite{smart}, and P.~Panzone, L.~Piovan, and M.~Ferrari \cite{panzone}.
Guinand \cite{guinand1}, \cite{guinand2} did state \eqref{2.11}.  Due to a
miscalculation, $\zeta(2n+1)$ did not appear in  Apostol's formula
\cite{apostol1}, but he later \cite{apostol2} realized his mistake and so
discovered \eqref{2.11}.  Also, recall that in Section \ref{sect2}, we
mentioned the very general transformation formula for multiple Barnes zeta
functions by Komori, Matsumoto, and Tsumura \cite{kmt} that contains Ramanujan's formula \eqref{2.11}
for $\zeta(2n+1)$  as a special case.  Lastly, in the beautiful
work of M.~Katsurada \cite{mk}, \cite{mk2} on asymptotic expansions of
$q$-series, Ramanujan's formula \eqref{i3.21} arises as a special case.

We also have not considered further
formulas for $\zeta(2n+1)$ that would arise from the differentiation of,
e.g., \eqref{2.6} and \eqref{2.10} with respect to $z$; see \cite{rocky}.

\section{The Associated Polynomials and their Roots}\label{sec:roots}

In this section we discuss the polynomials that are featured in Ramanujan's
formula \eqref{i3.21} and discuss several natural generalizations. These
polynomials have received considerable attention in the recent literature. In
particular, several papers focus on the location of zeros of these
polynomials. We discuss some of the developments starting with Ramanujan's
original formula and indicate avenues for future work.

Following \cite{gunmurtyrath} and \cite{murtysmythwang}, define the
{\emph{Ramanujan polynomials}}
\begin{equation}
  R_{2 m + 1} (z) = \sum_{n = 0}^{m + 1} \frac{B_{2 n}}{(2 n) !}  \frac{B_{2
  m - 2 n + 2}}{(2 m - 2 n + 2) !} z^{2 n} \label{eq:R}
\end{equation}
to be the polynomials appearing on the right-hand side of Ramanujan's formula
\eqref{i3.21} and \eqref{eq:rama:F}. The discussion in
Section~\ref{sec:eichler} demonstrates that the Ramanujan polynomials are,
essentially, the odd parts of the period polynomials attached to Eisenstein
series of level~$1$. In \cite{murtysmythwang}, M.R.~Murty, C.J.~Smyth and
R.J.~Wang prove the following result on the zeros of the Ramanujan
polynomials.

\begin{theorem}
  \label{thm:RR:roots}For $m \geq 0$, all nonreal zeros of the Ramanujan
  polynomials $R_{2 m + 1} (z)$ lie on the unit circle.
\end{theorem}

Moreover, it is shown in \cite{murtysmythwang} that, for $m \geq 1$, the
polynomial $R_{2 m + 1} (z)$ has exactly four real roots and that these
approach $\pm 2^{\pm 1}$, as $m \to\infty$.

\begin{example}
  As described in \cite{murtysmythwang}, one interesting consequence of this
  result is that there exists an algebraic number $\alpha \in \mathcal{H}$
  with $| \alpha | = 1$, $\alpha^{2 m} \neq 1$, such that $R_{2 m + 1} (\alpha) = 0$ and, consequently,
  \begin{equation}
    \frac{\zeta (2 m + 1)}{2} = \frac{F_{2 m + 1} (\alpha) - \alpha^{2 m}
    F_{2 m + 1} \left(- \frac{1}{\alpha} \right)}{\alpha^{2 m} - 1},
    \label{eq:zeta:Fa}
  \end{equation}
  where $F$ is as in \eqref{eq:Fa}. In other words, the odd zeta values can be
  expressed as the difference of two Eichler integrals evaluated at special
  algebraic values. An extension of this observation to Dirichlet $L$-series
  is discussed in \cite{bs}. Equation \eqref{eq:zeta:Fa} follows
  directly from Ramanujan's identity, with $z = \alpha$, in the form
  \eqref{eq:rama:F}. Remarkably, equation \eqref{eq:zeta:Fa} gets close to
  making a statement on the transcendental nature of odd zeta values: S.~Gun,
  M.R.~Murty and P.~Rath prove in \cite{gunmurtyrath} that, as
  $\alpha$ ranges over {\emph{all}} algebraic values in $\mathcal{H}$ with
  $\alpha^{2 m} \neq 1$, the set of values obtained from the right-hand side
  of \eqref{eq:zeta:Fa} contains at most one algebraic number.
\end{example}

As indicated in Section~\ref{sec:eichler}, the Ramanujan polynomials are the
odd parts of the period polynomials attached to Eisenstein series (or, rather,
period functions; see Remark~\ref{rk:pp} and \cite{zagier91}). On the
other hand, it was conjectured in \cite{lalin1}, and proved in
\cite{lalin2}, that the full period polynomial is in fact
{\emph{unimodular}}, that is, all of its zeros lie on the unit circle.

\begin{theorem}
  \label{thm:RRf:roots}For $m > 0$, all zeros of the polynomials
  \begin{equation*}
    R_{2 m + 1} (z) + \frac{\zeta (2 m + 1)}{(2 \pi i)^{2 m + 1}} (z^{2 m
     + 1} - z)
  \end{equation*}
  lie on the unit circle.
\end{theorem}

An analog of Theorem~\ref{thm:RR:roots} for cusp forms has been proved by
J.B.~Conrey, D.W.~Farmer and {\"O}.~Imamoglu \cite{cfi}, who show that,
for any Hecke cusp form of level $1$, the odd part of its period polynomial
has trivial zeros at $0$, $\pm 2^{\pm 1}$ and all remaining zeros lie on the
unit circle. Similarly, A.~El-Guindy and W.~Raji \cite{egr} extend
Theorem~\ref{thm:RRf:roots} to cusp forms by showing that the full period
function of any Hecke eigenform of level $1$ has all its zeros on the unit
circle.

\begin{remark}
  In light of \eqref{eq:pp:G} it is also natural to ask whether the
  polynomials
  \begin{equation*}
    p_m (z) = \frac{\zeta (2 m + 1)}{2} (1 - z^{2 m}) - \frac{(2 \pi i)^{2
     m + 1}}{2} \sum_{n = 1}^m \frac{B_{2 n}}{(2 n) !}  \frac{B_{2 m - 2 n +
     2}}{(2 m - 2 n + 2) !} z^{2 n - 1}
  \end{equation*}
  are unimodular. Numerical evidence suggests that these period polynomials
  $p_m (z)$ are indeed unimodular. Moreover, let $p_m^- (z)$ be the odd part
  of $p_m (z)$. Then, the polynomials $p_m^- (z) / z$ also appear to be
  unimodular. It seems reasonable to expect that these claims can be proved
  using the techniques used for the corresponding results in
  \cite{murtysmythwang} and \cite{lalin2}. We leave this question for the
  interested reader to pursue.
\end{remark}

Very recently, S.~Jin, W.~Ma, K.~Ono and K.~Soundararajan
\cite{jmos} established the following extension of the result of
El-Guindy and Raji to cusp forms of any level.

\begin{theorem}
  For any newform $f \in S_k (\Gamma_0 (N))$ of even weight $k$ and level
  $N$, all zeros of the period polynomial, given by the right-hand side of
  \eqref{eq:rw}, lie on the circle $| z | = 1 / \sqrt{N}$.
\end{theorem}

Extensions of Ramanujan's identity, and some of its ramifications, to higher
level are considered in \cite{bs}. Numerical evidence suggests
that certain polynomials arising as period polynomials of Eisenstein series
again have all their roots on the unit circle. Especially in light of the
recent advance made in \cite{jmos}, it would be interesting to study
period polynomials of Eisenstein series of any level more systematically.
Here, we only cite one conjecture from \cite{bs}, which concerns
certain special period polynomials, conveniently rescaled, and suggests that
these are all unimodular.

\begin{conjecture}
  \label{conj:unimod}For nonprincipal real Dirichlet characters $\chi$ and
  $\psi$, the {\emph{generalized Ramanujan polynomial}}
  \begin{equation}
    R_k (z ; \chi, \psi) = \sum_{s = 0}^k \frac{B_{s, \chi}}{s!} \frac{B_{k -
    s, \psi}}{(k - s) !} \left(\frac{z - 1}{M} \right)^{k - s - 1} (1 - z^{s
    - 1}) \label{eq:rx}
  \end{equation}
  is unimodular, that is, all its roots lie on the unit circle.
\end{conjecture}

Here, $B_{n, \chi}$ are the generalized Bernoulli numbers, which are defined
by
\begin{equation}
  \sum_{n = 0}^{\infty} B_{n, \chi}  \frac{x^n}{n!} = \sum_{a = 1}^L
  \frac{\chi (a) x e^{a x}}{e^{L x} - 1}, \label{eq:Bchi}
\end{equation}
if $\chi$ is a Dirichlet character modulo $L$. If $\chi$ and $\psi$ are both
nonprincipal, then $R_k (z ; \chi, \psi)$ is indeed a polynomial. On the other
hand, as shown in \cite{bs}, if $k > 1$, then $R_{2 k} (z ; 1, 1)
= R_{2 k - 1} (z) / z$, that is, despite their different appearance, the
generalized Ramanujan polynomials reduce to the Ramanujan polynomials,
introduced in \eqref{eq:R}, when $\chi = 1$ and $\psi = 1$. For more details
about this conjecture, we refer to \cite{bs}.

\section{Further History of Proofs of \eqref{i3.21}}\label{sect4}

In this concluding section, we shall offer further sources where additional
proofs of Ramanujan's formula for $\zeta(2n+1)$ can be found.  We emphasize
that two papers of Grosswald \cite{gross1}, \cite{gross2}, each proving
formulas from Ramanujan's second notebook, with the first giving a proof of
\eqref{i3.21}, stimulated the first author and many others to seriously
examine the content of Ramanujan's notebooks \cite{nb}. K.~Katayama
\cite{katayama1}, \cite{katayama1}, H.~Riesel \cite{riesel}, S.N.~Rao
\cite{snrao},  L.~Vep\u{s}tas \cite{vepstas}, N.~Zhang \cite{zhang}, and
N.~Zhang and S.~Zhang \cite{zhangzhang} have also developed proofs of
\eqref{i3.21}.  B.~Ghusayni \cite{ghusayni} has used Ramanujan's formula for
$\zeta(2n+1)$ for numerical calculations.

For further lengthy  discussions, see the first author's book
\cite[p.~276]{II}, his  survey paper \cite{survey}, his account of
Ramanujan's unpublished manuscript \cite{jrms}, and his fourth book with
Andrews \cite[Chapter 12]{geabcbIV} on Ramanujan's lost notebook.  These
sources also contain a plethora 
of references to proofs of the third and fourth companions to Ramanujan's
formula for $\zeta(2n+1)$.  Another survey has been given by S.~Kanemitsu and
T.~Kuzumaki \cite{kk}.

Also, there exist a huge number of generalizations and analogues of
Ramanujan's formula for $\zeta(2n+1)$.  Some of these are discussed in
Berndt's book \cite[p.~276]{II} and his paper \cite{crelle}.  S.-G.~Lim
\cite{lim1}, \cite{lim2}, \cite{lim3} has established an enormous number of
beautiful identities in the spirit of the identities discussed in the present
survey.  Kanemitsu, Y.~Tanigawa, and M.~Yoshimoto \cite{kty}, \cite{kty2},
interpreting Ramanujan's formula \eqref{i3.21} as a modular transformation,
derived further formulas for $\zeta(2n+1)$, which they have shown lead to a
rapid calculation of $\zeta(3)$ and $\zeta(5)$, for example.

In Sections~\ref{sec:eisenstein} and \ref{sec:eichler} we have seen that
Ramanujan's formula can be viewed as describing the transformation laws of the
modular Eisenstein series $E_{2 k}$, where $k > 1$, of level $1$ (that is,
with respect to the full modular group), the quasimodular Eisenstein series
$E_2$ as well as their Eichler integrals. We refer to \cite{bs}
for extensions of these results to higher level.

\end{document}